\documentclass[12pt]{article}
\usepackage{latexsym}
\usepackage{amssymb}

\newtheorem{teo}{THEOREM}[section]
\newtheorem{prop}[teo]{PROPOSITION}

\newtheorem{obser}[teo]{Remark}
\newtheorem{defi}[teo]{DEFINITION}
\newtheorem{coro}[teo]{COROLLARY}
\newenvironment{proof}{\noindent \sf Proof. \sf}

\def\Hom{\mathop{\rm Hom}\nolimits}

\def\qed{\hfill \mbox{$\square$}\bigskip}

\def\Ext{\mathop{\rm Ext}\nolimits}
\def\Tor{\mathop{\rm Tor}\nolimits}

\def\Ker{\mathop{\rm Ker}\nolimits}

\def \Z{{\mathbb Z}}
\def \Q {{\mathbb Q}}
\def\dim{\mathop{\rm dim_k}\nolimits}
\def\dimQ{\mathop{\rm dim_\Q}\nolimits}
\def\rank{\mathop{\rm rank}\nolimits}

\def\ootimes{\mbox{\footnotesize$\otimes$}}

\begin{document}
\sf
\title{Cartan--Leray spectral sequence for Galois coverings of linear categories}
\author{Claude Cibils and Mar\'\i a Julia Redondo \footnote{The authors gratefully acknowledge PICS 1514 of the CNRS. The second author had a position of
\emph{chercheur associ\'e} from the CNRS at Montpellier during the preparation of this work.}}
\date{}

\maketitle

\begin{abstract}
We provide a Cartan-Leray type spectral sequence for the
Hochschild-Mitchell (co)homology of a Galois covering of linear
categories. We infer results relating the Galois group and
Hochschild cohomology in degree one.
\end{abstract}

\small \noindent 2000 Mathematics Subject Classification : 16E40, 18E05

\noindent Keywords : Cartan, Leray, covering, Hochschild, Mitchell, category, homology, cohomology.

\font\fivrm=cmr5 \relax
\input{prepictex}
\input{pictex}
\input{postpictex}

\section {\sf Introduction }

This past years several results and tools from algebraic topology
has been adapted to representation theory of non commutative
associative finite dimensional algebras over a field $k$. For this
purpose $k$-categories has been considered as Galois coverings of
algebras by P. Gabriel, see for instance \cite{ga} or also
\cite{bo-ga}. Recall that a $k$-category is a small category where
morphisms are $k$-vector spaces and composition is $k$-bilinear.
Relations between the representation theory of a $k$-algebra and
functors from the universal cover has also been described and the
fundamental group of the presentation of an algebra by a quiver
with relations emerged in this context, see \emph{eg.} \cite{ga}.

Moreover a precise relation which can be compared with Hurewicz's Theorem, (see for instance \cite{m-c}), has been
obtained in \cite{a-jap} and \cite{jap-s} between the fundamental group of a presentation and Hochschild cohomology
of the algebra in degree one. More recently an interpretation of this algebraic fundamental group as a fundamental
group of a simplicial complex or of a classifying CW-complex has been described, see \cite{bu1,bu2,re}. Note also
that comparisons between Hochschild cohomologies of a Galois covering are obtained in \cite{jap-m}. In \cite{re-ri}
the relation between skew group algebras and coverings is studied and results are provided for a finite Galois group
in semi-simple characteristic.

Our main purpose in this paper is to provide a setting relating
the Galois group of the covering of a $k$-algebra (or more
generally of a $k$-category) with the Hochschild (co)homology
theory of the situation. A model for this is the Cartan-Leray
spectral sequence associated to a group acting freely and properly
on a connected space, see \cite{c-l,c} and for instance \cite[p.
337]{m-c} or \cite[p. 206]{w}.

We provide a Cartan-Leray type spectral sequence for Galois
coverings of $k$-categories. For this purpose we briefly recall
the so-called Hochschild-Mitchell homology and cohomology theories
$H_*(\cal C, M)$ and $H^*(\cal C, M)$ of a $k$-category $\cal C$
with coefficients in a bimodule $\cal M$ over $\cal C$. These
(co)homology theories have been introduced by Mitchell in
\cite{m}, and are closely related with theories considered in
\cite{ke,mcca} but differs from Quillen's homology described in
\cite[Appendix C]{lo}. Notice that the spectral sequence results
we obtain are also valid in case of a commutative ring $k$ instead
of a field, provided that we consider $k$-categories with flat
module morphisms. Hochschild-Mitchell (co)homology theories are
used in theoretical informatics and we thank speakers and
participants of a seminar in Montpellier related to the subject
for pointing out this theory.

In case of a $k$-category with a finite number of objects the
Hochschild-Mitchell (co)homology coincides with the usual
Hochschild (co)homology of the $k$-algebra associated to the
category, namely the direct sum of the morphisms spaces equipped
with the matrix product. We provide a proof of this well known
agreement property, which in turn enables us to prove that
Hochschild-Mitchell theory is invariant under contraction and
expansion of a $k$-category - precise definitions of these
operations are provided in the text.

A difference with the Cartan-Leray spectral sequence in the algebraic topology context is the use of coefficients.
For Hochschild-Mitchell (co)homology, coefficients are taken in a bimodule which has to be lifted through the
canonical projection, while in the topological context an unchanged abelian group of coefficients is used.

A cohomological type spectral sequence cannot be obtained in
general for an infinite category with a group acting freely.
Indeed Hochschild-Mitchell cohomology makes use of products of
infinitely many $k$-nerves of the category and we note that this
cannot be avoided. The action of the group on those cochains do
not provide free modules over the group algebra, hence the
spectral sequence does not necessarily collapse for the column
filtration.

Nevertheless our interest is mainly in Hochschild cohomology of a
finite dimensional algebra with coefficients in itself. An usual
duality holds between Hochschild-Mitchell homology and cohomology
for locally finite bimodules, that is, bimodules $\cal{M}$ with
values in the category of finite dimensional vector spaces. The
dual $D\cal M$ of any bimodule is well defined and in case $\cal
M$ is locally finite we have $DD\cal M = \cal M$.

This duality provides a Cartan-Leray type spectral sequence of
cohomological type with coefficients in a locally finite bimodule.
As a particular case of interest we consider a finite dimensional
$k$-algebra $\Lambda$ such that the larger semi-simple quotient is
a product of copies of $k$, provided with a complete system of
primitive orthogonal idempotents. Let $\cal {B}$ be the
corresponding $k$-category and let $\cal C \to \cal {B}$ be a
Galois covering with group $G$. Then the spectral sequence
provides a canonical embedding $\Hom(G, k^+)$ in $H^1(\cal B, \cal
B)$. Note that this result has been obtained by Assem and de la
Pe\~na in \cite{a-jap} in case $G$ is the fundamental group of a
presentation of $\cal B$ by a triangular quiver with admissible
relations, afterwards de la Pe\~na and Saor\'\i n noticed in
\cite{jap-s} that the triangular assumption is not needed.

As a consequence of our result the dimension of $H^1(\cal B, \cal
B)$ is bounded below by the rank of $G$, compare with \cite{jap-m}
where this inequality is obtained in case of a free group $G$
acting freely on a schurian category $\cal C$.

In \cite{jap-m} a comparison between the Hochschild cohomology of
the algebras involved in a Galois covering $A \to B$ with a cyclic
group of order $2$ is initiated, where $B$ is the Kronecker
algebra. Our considerations show that for a field of
characteristic not $2$ the $3$--dimensional vector space $H^1(B)$
is actually isomorphic to the fixed points of the $5$--dimensional
vector space $H^1(A,LB)$, where $LB$ is the lifted bimodule.

More generally let $\cal C \to \cal B$ be a Galois covering of
finite categories with finite group $G$, and assume that the
characteristic of $k$ is zero or do not divide the order of $G$.
We infer from the Cartan-Leray spectral sequence that
$$H^n({\cal B}, {\cal M }) = H^n ({\cal C}, L {\cal M})^G$$
where $L \cal M$ is the lifting of the $\cal B$-bimodule $\cal M$
to the category $\cal C$.

\section{\sf Hochschild--Mitchell (co)homology}\label{sectionHOMI}

A category is called \emph{small} if the objects and morphisms are
sets, all categories in this paper are small. Let $k$ be a field,
a \emph{$k$-category} is a small category such that the morphisms
are $k$-vector spaces and the composition of morphisms is
$k$-bilinear.

Let $\cal C$ be a $k$-category with objects set ${\cal C}_0$. It
is convenient to denote by ${ }_y{\cal C}_x$ the morphisms from
the object $x$ to the object $y$, composition provides $k$-linear
maps
$${ }_z{\cal C}_y \otimes_k { }_y{\cal C}_x \to { }_z{\cal C}_x$$
for $x,y,z \in {\cal C}_0$.

A \emph{bimodule} $\cal M$ over a $k$-category $\cal C$ is a
bifunctor from ${\cal C} \times {\cal C}^{op}$ to the category of
vector spaces. In other words $\cal M$ is given by a set of vector
spaces $\{ { }_y{\cal M}_x \}_{x,y \in {\cal C}_0}$ and left and
right actions

\[ \begin{array}{ccccccc}
{ }_z{\cal C}_y & \otimes & { }_y{\cal M}_x & & &\to &{ }_z{\cal M}_x\\
& & { }_y{\cal M}_x & \otimes & { }_x{\cal C}_u & \to & { }_y{\cal
M}_u
\end{array} \]
satisfying the usual associativity conditions, namely
\begin{eqnarray*}
{ }_uc_z \ ({ }_zc_y \ { }_ym_x) & = & ({ }_uc_z\ { }_zc_y)\ { }_ym_x \\
{ }_ym_x \ ({ }_xc_u \ { }_uc_v) & = & ({ }_ym_x\ { }_xc_u)\ { }_uc_v\\
({ }_zc_y \ { }_ym_x ) \ { }_xc_u & = & { }_zc_y \ ({ }_ym_x \ { }_xc_u)\\
{ }_y1_y \ { }_ym_x & = & { }_ym_x \ { }_x1_x \ = \ { }_ym_x.
\end{eqnarray*}
For instance the standard bimodule over a $k$-category is the
category itself.

Of course if $\cal C$ is finite the vector space
$$\Lambda_{\cal C} = \bigoplus_{x,y \in {\cal C}_0} { }_y{\cal C}_x$$
has an algebra structure, with a well defined matrix product given
by the composition of $\cal C$, the identity is the matrix with
${}_x1_x$ in the diagonal and zero otherwise. Note that
$\Lambda_{\cal C}$ has by construction a complete set of
orthogonal idempotents, namely $\{ { }_x1_x \}_{x \in {\cal C}}$.
Usual $\Lambda_{\cal C}$-bimodules and $\cal C$-bimodules
coincide.

Conversely, let $\Lambda$ be a $k$-algebra equipped with a
complete set $F$ of orthogonal idempotents (non necessarily
primitive), that is, a finite subset $F$ of $\Lambda$ whose
elements satisfy $x^2=x$, $xy=0$ if $x \not=y$ and $\sum_{x \in F}
x=1$. The associated finite category ${\cal C}_{\Lambda, F}$ has
set of objects $F$ and $y \Lambda x$ are the morphisms from $x$ to
$y$. Composition is given by the product in $\Lambda$.

\begin{defi}
Let $(x_{n+1}, x_n, \dots , x_1)$ be a $n+1$-sequence of objects
of a $k$-category $\cal C$. The $k$-nerve associated to this
sequence is the vector space
$${ }_{x_{n+1}}{\cal C}_{x_n} \otimes \dots \otimes { }_{x_2}{\cal C}_{x_1}.$$
The $k$-nerve ${\cal N}_n$ of degree $n$ is the direct sum of all
the $k$-nerves associated to $n+1$-sequences of objects. We have
\begin{eqnarray*}
{\cal N}_0 & = & \bigoplus_{x \in {\cal C}_0} { }_x {\cal C}_x \\
{\cal N}_n & = & \bigoplus_{ \mbox{\tiny$n+1$-sequences}} {
}_{x_{n+1}}{\cal C}_{x_n} \otimes \dots \otimes { }_{x_2}{\cal
C}_{x_1}.
\end{eqnarray*}
\end{defi}

\begin{defi}
Let $\cal C$ be a $k$-category and let $\cal M$ be a $\cal
C$-bimodule. The Hochschild-Mitchell cohomology $H^*({\cal
C},{\cal M})$ is the cohomology of the cochain complex
$$ 0 \longrightarrow \prod_{x} { }_x{\cal M}_x \stackrel{d}{\longrightarrow}
\Hom ({\cal N}_1 , {\cal M}) \stackrel{d}{\longrightarrow} \dots
\stackrel{d}{\longrightarrow} \Hom ({\cal N}_n, {\cal M})
\stackrel{d}{\longrightarrow} \dots $$ where
$$C^n( {\cal C}, {\cal M}) = \Hom ({\cal N}_n, {\cal M}) = \prod_{ \mbox{\tiny $n+1$-sequences}}
\Hom_k ( { }_{x_{n+1}}{\cal C}_{x_n} \otimes \dots \otimes {
}_{x_2}{\cal C}_{x_1},\ { }_{x_{n+1}}{\cal M}_{x_1})$$ and $d$ is
given by the usual formulas of Hochschild cohomology as follows.
Let $f=\left\{f_{\left(x_{n+1},\dots,\
x_1\right)}\right\}_{\mbox{\tiny $n+1$-sequences}}$ be a family of
linear maps. Then $df$ is the family of linear maps

$$(df)=\left\{df_{{\left(x_{n+2},\cdots,\
x_1\right)}}\right\}_{\mbox{\tiny $n+2$-sequences}}\mbox { where
}$$
$$\begin{array}{l} df_{\left(x_{n+2},\dots,\
x_1\right)}\left({}_{x_{n+2}}c_{x_{n+1}}\ootimes\cdots\ootimes{}_{x_2}c_{x_1}\right)=\\
(-1)^{n+1}\left({}_{x_{n+2}}c_{x_{n+1}}\right)f_{\left(x_{n+1},\dots,\
x_1\right)}
\left({}_{x_{n+1}}c_{x_{n}}\ootimes\cdots\ootimes{}_{x_2}c_{x_1}\right)+\\
\sum_{i=1}^{i=n}(-1)^{i}f_{\left(x_{n+2},\dots,\ x_{i+2},\
x_{i},\dots\ x_1\right)}
\left({}_{x_{n+2}}c_{x_{n}}\ootimes\cdots\ootimes\left({}_{x_{i+2}}c_{x_{i+1}}\right)\left({}_{x_{i+1}}c_{x_{i}}\right)
\ootimes\cdots\ootimes{}_{x_2}c_{x_1}\right)+\\
f_{\left(x_{n+2},\dots,\ x_2\right)}
\left({}_{x_{n+2}}c_{x_{n+1}}\ootimes\cdots\ootimes{}_{x_3}c_{x_2}\right)\left({}_{x_{2}}c_{x_{1}}\right).
\end{array}$$

Note that each linear map of the family $df$ is well defined.
\end{defi}

\begin{obser}\label{centro}
The tensor product of $\cal C$-bimodules ${\cal M}_1$ and ${\cal
M}_2$ is the bimodule given by
$${ }_y({\cal M}_1 \otimes {\cal M}_2)_x = \bigoplus_{z} { }_y({\cal M}_1)_z \otimes
{ }_z({\cal M}_2)_x .$$

As in the algebra case, the standard bimodule $\cal C$ has a
resolution by tensor powers of $\cal C$ which are projective
bimodules. Applying the functor $\Hom (-, {\cal M})$ provides
precisely the cochain complex above. Consequently
Hochschild-Mit\-chell cohomology is an instance of an $\Ext$
functor. In particular
$$H^0({\cal C}, {\cal M}) = \{ ({ }_xm_x)_x \ \mid
\mbox{ ${ }_xm_x \in { }_x{\cal M}_x$ and ${ }_yf_x \ { }_xm_x = {
}_ym_y \ {}_yf_x$ for all ${ }_yf_x \in {}_y{\cal C}_x $} \}.$$
Hence $H^0({\cal C}, {\cal C})$ is the center of the category.
\end{obser}

\begin{obser}
Replacing the direct product by a direct sum in the cochain
complex has no sense since the Hochschild coboundary will provide
in general a direct product of morphisms, even as an image of a
morphism located in a single component. However if the category is
locally finite the direct sums of morphisms at the cochain level
provide a subcomplex.
\end{obser}

\begin{defi}
Let $\cal C$ be a $k$-category and let $\cal M$ be a $\cal
C$-bimodule. The Hochschild-Mitchell homology $H_*({\cal C},{\cal
M})$ is the homology of the chain complex
$$ \dots \longrightarrow {\cal M} \otimes {\cal N}_n \stackrel{d}{\longrightarrow}
\dots \stackrel{d}{\longrightarrow} {\cal M} \otimes {\cal N}_1
\longrightarrow \bigoplus_{x} {}_x{\cal M}_x \longrightarrow 0 $$
where
$$C_n({\cal C},{\cal M}) = {\cal M} \otimes {\cal N}_n =
\bigoplus_{\mbox{\tiny $n\!+\!1$-sequences}}{ }_{x_1}{\cal
M}_{x_{n+1}} \otimes { }_{x_{n+1}}{\cal C}_{x_n} \otimes \dots \otimes
{ }_{x_2}{\cal C}_{x_1}$$ and $d$ is given by the usual Hochschild
boundary, see for instance \cite{lo}.
\end{defi}

Again this is $\Tor_*({\cal C},{\cal M})$ in the abelian
category of $\cal C$-bimodules.

\begin{obser}
Hochschild-Mitchell homology and cohomology theories are not
homology or cohomology of small categories as considered in
\cite[Appendix C]{lo}. The following well known agreement result
betwen usual Hochschild and Hochschild-Mitchell (co)homology
theory holds.
\end{obser}

\begin{prop}
Let $\cal C$ be a finite $k$-category, $\Lambda_{\cal C}$ the
corresponding algebra and ${\cal M}$ a $\Lambda_{\cal
C}$-bimodule. Then
$$H^*({\cal C}, {\cal M}) = H^*(\Lambda_{\cal C}, {\cal M}) \quad \mbox{and} \quad
H_*({\cal C}, {\cal M}) = H_*(\Lambda_{\cal C}, {\cal M})$$ where
$H(\Lambda_{\cal C}, {\cal M})$ denotes usual Hochschild
(co)homology of the algebra $\Lambda_{\cal C}$.

\end{prop}

\begin{proof}
Let $\Lambda := \Lambda_{\cal C}$ and consider the semi-simple
subalgebra of $\Lambda$ given by $E= \prod k \ {}_x1_x$. Note that
any $E$-bimodule is projective since the enveloping algebra of $E$
is still semi-simple. Hence $\Lambda \otimes_E X \otimes_E
\Lambda$ is a projective $\Lambda$-bimodule for any $E$-bimodule
$X$. Consequently there is a projective resolution of $\Lambda$ as
a $\Lambda$-bimodule given by
$$ \dots \longrightarrow \Lambda \otimes_E \Lambda \otimes_E \dots \otimes_E \Lambda
\longrightarrow \dots \longrightarrow \Lambda \otimes_E \Lambda
\longrightarrow \Lambda \longrightarrow 0.$$ Applying the functor
$\Hom_{\Lambda-\Lambda}(-, {\cal M})$ to this resolution and
considering the canonical vector-space isomorphism
$$\Hom_{\Lambda-\Lambda}(\Lambda \otimes_E X \otimes_E \Lambda, {\cal M}) =
\Hom_{E-E}(X, {\cal M})$$ we obtain a cochain complex computing
$H^*(\Lambda, {\cal M})$ which coincides with the complex we have
defined for the Hochschild-Mitchell cohomology. The same argument
shows the assertion for homology. \qed
\end{proof}

This result has an immediate generalization that we provide
although we will not make direct use of it, however it can serve
as a guideline for the reader. Notice first that there are two
operations that we can perform on a $k$-category:

Let $\cal F$ be a finite full subcategory of a $k$-category $\cal
C$. The \emph{contraction} of $\cal C$ along $\cal F$ is the
category ${\cal C}_{\cal F}$ where $\cal F$ is replaced with one
object having $\Lambda_{\cal F}$ as endomorphism algebra and the
evident morphisms between this new object and the rest of the
objects of $\cal C$.

The \emph{expansion} of $\cal C$ along a complete system $F$ of
orthogonal idempotents of ${}_x {\cal C}_x$ for an object $x$ is
given by the category obtained by replacing the object $x$ by the
set $F$ and providing morphisms in accordance.

Of course a bimodule $\cal M$ over $\cal C$ provides a uniquely
determined bimodule over a contraction or an expansion of $\cal C$
and clearly we have the following result.

\begin{prop}
The Hochschild-Mitchell (co)homology is invariant under
contraction and expansion.
\end{prop}

Finally we note that the cyclic Mitchell homology of a category is
also well defined along this lines.

\section {\sf Cartan-Leray spectral sequence}

Let $\cal C$ be a $k$-category and $G$ be a group acting freely on
it, namely there is an action on objects such that for any $s \in
G$, $x \in {\cal C}_0$, $s x = x $ implies $s=1$, and an action on
the morphisms such that for any $c \in {}_y{\cal C}_x$ we have $s
c \in{}_{sy}{\cal C}_{sx}$ and $s (dc) = sd \ sc$ for any $d \in
{}_z{\cal C}_y$, $c \in {}_y{\cal C}_x$. Of course $s(tx)=(st)x$
and $s(tc)=(st)c$ for any $s,t \in G$, $x \in {\cal C}_0$ and $c$
a morphism in $\cal C$.

The quotient category ${\cal C}/G$ has objects the set of orbits
of objects. The morphisms between two orbits is the direct sum of
the orbit spaces of morphisms, namely let $u,v$ be orbits of
objects, then $\left( \bigoplus_{x \in u,\ y \in v} { }_y{\cal
C}_x\right)$ is a left $kG$-module and we put
$${}_v ({\cal C}/G)_u = \left( \bigoplus_{x \in u,\ y \in v} { }_y{\cal C}_x\right)/G$$
where as usual we denote $X/G$ the quotient of a $kG$-module $X$
by the action of $G$, namely $X/G = X/ (\Ker \epsilon) X$ where
$\epsilon:kG \to k$ is the augmentation algebra morphism defined
by $\epsilon(s)=1$ for all $s \in G$. In particular $kG/G=k$.

Composition of morphisms is well defined precisely because the
$G$-action on objects of $\cal C$ is free.

\begin{obser}
P. Gabriel has considered in \cite[p. 85]{ga} a group $G$ acting
freely on a $k$-category and fixed points on families of
morphisms. The fixed points approach impose restrictions on the
categories considered in \cite{ga} which we do not need in our
context. In the context of bocses  Y.A. Drozd and S.A. Ovsienko
\cite{drov} has introduced the same categorical quotient that we
consider in this paper, we thank L. Salmer\'on for pointing out
this.

\end{obser}

\begin{obser}\label{galois}
Let $x_0$ be a fixed object in an orbit $u$, in other words
$u=\overline{x_0}$ where $\overline{x_0}$ denotes the $G$-orbit of
an object of $\cal C$. Then
$${}_v ({\cal C}/G)_u = \bigoplus_{y \in v} { }_y{\cal C}_{x_0}.$$
\end{obser}

\begin{defi}
A Galois covering of $k$-categories under the action of a group
$G$ is a functor ${\cal C} \to {\cal B}$ where $G$ acts freely on
$\cal C$ and ${\cal B} = {\cal C}/G$. The functor is the
projection functor.
\end{defi}

\begin{defi}
Let $\cal M$ be a $\cal B$-bimodule. The lifted $\cal C$-bimodule
$L{\cal M}$ is defined by
$${ }_y(L{\cal M})_x = { }_{\overline{y}}{\cal M}_{\overline{x}}$$
with left action
$${ }_z {\cal C}_y \otimes { }_y(L{\cal M})_x \to { }_z(L{\cal M})_x$$
given by ${ }_zc_y \ { }_ym_x \ = \ { }_z({
}_{\overline{z}}\overline{c}_{\overline{y}} \ {
}_{\overline{y}}m_{\overline{x}})_x$, and right action defined
analogously.
\end{defi}

Note that this definition is in force for any functor between
categories and, in particular, for a Galois covering. However for
a Galois covering the group $G$ acts on a $\cal C$-bimodule $\cal
U$: the $\cal C$-bimodule $ {}^s \cal U$ is given by the vector
spaces $ {}_y\left({}^s \cal U\right)_x = {}_{s^{-1}y} {\cal
U}_{s^{-1}x}$ and left action
$$c . u = \left(s^{-1}c\right)(u) \in {}_{s^{-1}z} {\cal U}_{s^{-1}x},
\quad c \in  {}_z{\cal C}_y,$$
and similarly on the right. Clearly a lifted bimodule is fixed
under this action.

\begin{prop}
Let $\cal U$ be a $\cal C$-bimodule. Then for each $s\in G$ there
is an induced map between Hochschild-Mitchell homologies
$$H_{*}\left(\cal C, \ \cal U\right) \longrightarrow H_{*}\left({\cal C}, \ {}^{s}\cal U\right)$$
\end{prop}

\begin{proof}
Consider at the Hochschild-Mitchell chain level the map given as
follows, (tensor signs are replaced by commas)
$$\left({ }_{x_{1}}{(u)}_{x_{n+1}},\ { }_{x_{n+1}}({c_n})_{x_{n}},
\dots, { }_{x_{2}}({c_1})_{x_{1}}\right)\mapsto \left(u,\
sc_n,\dots, sc_1\right),$$ where $u$ on the right hand side
belongs to ${}_{sx_{1}}({}^{s} {\cal U})_{sx_{n+1}}=
{}_{x_{1}}\left(\cal U\right)_{x_{n+1}}.$ This is clearly a chain
map as a direct consequence of the definitions of left end right
actions on ${}^{s}{\cal U}$.
\end{proof}

\begin{coro}
Let $\cal M$ be a $\cal B$-bimodule and $L\cal M$ the
corresponding lifted $\cal C$-bimodule. Then the Galois group $G$
of the covering acts on $H_*({\cal C}, L\cal M).$
\end{coro}

\begin{proof}
As already quoted we have ${}^s(L{\cal M})= L\cal M$ for any $s\in
G$.\qed
\end{proof}



Consequently $H_n({\cal C},L{\cal M})$ is a $kG$-module for
each $n \geq 0$.

\begin{prop}
Let ${\cal C} \to {\cal B}={\cal C}/G$ be a Galois covering of
$k$-categories with group $G$ and let $\cal M$ be a $\cal
B$-bimodule. Then $C_*({\cal C}, L{\cal M})$ is a free $kG$-module
and $$C_*({\cal C},L{\cal M})/G = C_*({\cal B}, {\cal M}).$$
\end{prop}
\begin{proof}
The $k$-nerve of degree $n$ of $\cal C$ is a free $kG$-module
which decomposes as a direct sum of $kG$-modules along $G$-orbits
of $n+1$ objects, each of the summands is a free $kG$-module. The
Hochschild-Mitchell chain complex corresponding to the orbit of a given
(n+1)-sequence $x_{n+1}, \cdots, x_1$ is obtained by tensoring the
summand of the nerve with the constant vector space ${}_{\overline
x_1}{\cal M}_{\overline x_{n+1}}$, providing in this way that each
$kG$-direct summand of $C_*({\cal C}, L{\cal M})$ is
$kG$-free.

Let now $u_{n+1},\cdots,u_{1}$ be a sequence of $G$-orbits of objects of
$\cal C$. We assert that the two following vector spaces are
isomorphic,
$$A = \left[\bigoplus_{x_i\in u_i}\left( {}_{x_1}(L{\cal M})_{x_{n+1}}\otimes
{}_{x_{n+1}}{\cal C}_{x_n}\otimes \cdots \otimes {}_{x_{3}}{\cal
C}_{x_2}\otimes {}_{x_{2}}{\cal C}_{x_1} \right)\right]/G$$ and
$$B= {}_{u_1}{\cal M}_{u_{n+1}}\otimes
{}_{u_{n+1}}({\cal C}/G)_{u_n}\otimes \cdots \otimes
{}_{u_{3}}({\cal C}/G)_{u_2}\otimes {}_{u_{2}}({\cal C}/G)_{u_1}.
$$

Let $(m,c_n,\cdots, c_1)$ be a representative of $A$, then
$$\varphi(m,c_n,\dots, c_1)=(m,{\overline c_n},\dots,{\overline
c_1})$$ provides a well defined element of $B$. Conversely let
$(m,\gamma_n,\cdots,\gamma_2,\gamma_1)\in B$ and we choose any
$c_1\in \gamma_1$. We adjust the choice of $c_2\in\gamma_2$ in
such a way that the source object of $c_2$ is the target one of
$c_1$, this can be done in an unique way since the action of $G$
on the category is free. We pursue the right choices until
$c_n\in\gamma_n$. Concerning $m\in{}_{u_1}{\cal M}_{u_{n+1}}$,
choose to consider it in ${}_{x_1}(L{\cal M})_{x_{n+1}}$ and put
$$\psi(m,\gamma_n,\cdots,\gamma_1)= (m,c_n,\dots,c_1)\in \left[{}_{x_1}(L{\cal M})_{x_{n+1}}\otimes
{}_{x_{n+1}}{\cal C}_{x_n}\otimes \dots \otimes {}_{x_{2}}{\cal
C}_{x_1}\right]/G.$$

The point is to prove that this element does not depend on the
choice we made of a representative $c_1$ of the class $\gamma_1$.
Let $c'_1$ be another choice. There exist a unique $s\in G$ such
that $c'_1=sc_1$. Clearly the unique $c'_2$ such that
its source is the target of $c'_1$ is given by $c'_2=sc_2$, the
same group element $s$ is used. Finally $c'_n=sc_n$ and $m$ is now
located in ${}_{sx_1}(L{\cal M})_{sx_{n+1}}$. Note that
$$(m,sc_n,\dots,sc_1)=s(m,c_n,\dots,c_1)$$ so the elements
coincide in $A$.\qed
\end{proof}
\begin{obser}
The proof above shows that it is crucial for the bimodule $L\cal
M$ to verify
$${}_{x_1}(L{\cal M})_{x_{n+1}}={}_{sx_1}(L{\cal M})_{sx_{n+1}}={}_{{\overline x_1}}
{\cal M}_{{\overline x_{n+1}}}.$$
\end{obser}

\begin{obser}
The dual statement for cohomology does not hold in general for an
infinite group. Namely Hochschild-Mitchell cochain $C^n({\cal C},
L{\cal M})$ is an infinite product of free $kG$-modules which is
not a free module in general.
\end{obser}

\begin{teo}\label{hoCL} [Cartan-Leray]
Let ${\cal C}\to {\cal B}$ be a Galois covering of $k$-categories
with group $G$ and let ${\cal M}$ be a $\cal B$-bimodule. There is
a first quadrant spectral sequence converging to $H_n({\cal
B},{\cal M})$ with level $2$ term
$$E_{p,q}^2 = H_p(G, H_q({\cal C}, L{\cal M}))$$
where $H_p(G,X)$ denotes usual group homology with coefficients in
a $kG$-module $X$.
\end{teo}

We will prove this Theorem following the classical proof in
algebraic topology. Before, we consider a Cartan-Leray Theorem for
Hochschild-Mitchell cohomology which is of special interest for
finite dimensional $k$-algebras.

\begin{teo}\label{cohoCL} [Cartan-Leray]
Let ${\cal C}\to {\cal B}$ be a Galois covering of $k$-categories
with group $G$ and let ${\cal M}$ be a locally finite $\cal
B$-bimodule. There is a first quadrant spectral sequence of
cohomological type converging to $H^n({\cal B},{\cal M})$ with
level $2$ term
$$E^{p,q}_2 = H^p(G, H^q({\cal C}, L{\cal M}))$$
where $H^p(G,X)$ denotes usual group cohomology with coefficients
in a $kG$-module $X$.
\end{teo}

\begin{proof}
For any $\cal B$-bimodule $\cal Y$ we have
$$H_n({\cal B, Y})^* = H^n({\cal B}, D{\cal Y})$$
where ${}^*$ denotes the dual vector space and $D\cal Y$ is the
dual bimodule of $\cal Y$. Recall that a localy finite bimodule
$\cal M$ has finite dimensional vector spaces ${}_y{\cal M}_x$ for
each couple of objects $(y,x)$. Hence $DD\cal M = M$.

The Cartan-Leray spectral sequence of homological type provides by
dualisation a spectral sequence of cohomological type as required,
noticing that the functors $D$ and $L$ commute.\qed

\end{proof}

\begin{proof} Theorem \ref{hoCL}.
Let $k$ be the trivial $kG$-module and let $P_{\cdot} \to k \to 0$
be any $kG$-projective resolution of $k$. Consider the double
complex whose $(p,q)$ spot is $P_p \otimes_{kG} C_q({\cal C},
L{\cal M})$ with the differentials given by the tensor product of
those of the resolution of $k$ and those of the
Hochschild-Mitchell complex.

The homology of each row is usual group homology of $G$ with
coefficients in a free $kG$-module. Hence the homology of the
$q$-row is concentrated in $p=0$ where we have
$$H_0(G, C_q({\cal C}, L{\cal M})) = C_q({\cal C}, L{\cal M})/G = C_q({\cal C}/G, {\cal M}).$$
The induced vertical differentials are the Hochschild-Mitchell
boundaries ones by construction. Hence the filtration by the rows
of the double complex provides at level $2$ a $p=0$ column of
vector spaces which are $H_q({\cal C}/G, {\cal M})$. The
differentials at this level come and go to zero, consequently the
abutment of the spectral sequence is the Hochschild-Mitchell
homology of ${\cal B}= {\cal C}/G$ with coefficients in $\cal M$
as required.

The homology of the $p$-column is $P_p \otimes_{kG}H_.({\cal C},
L{\cal M})$ with horizontal differentials given by the projective
resolution of $k$. In other words the homology of the rows at
level $1$ is usual homology of $G$, namely $H_p(G, H_q({\cal C},
L{\cal M}))$. \qed
\end{proof}

\begin{obser}
An alternative approach for obtaining the Cartan-Leray spectral
sequence of a Galois covering of $k$-categories is to use
Grothendieck spectral sequence associated to the composition of
two functors, see for instance \cite[p. 150]{w}, in which case one
have to check that the requirements of Grothendieck's Theorem are
fulfilled either in the context of this paper or in the algebraic
topology context. We thank Serge Bouc for pointing out this fact.
\end{obser}

\begin{obser}
Let $\cal C$ be a small category which is not necessarily a
$k$-category. The set theoretical nerve provides a simplicial set
and the Cartan-Leray spectral sequence we have obtained is
analogous to the spectral sequence of a fibration, see for
instance \cite[p. 157]{gz}. However we use in this paper
$k$-categories obtaining simplicial vector spaces instead of
simplicial sets. In other words the categories we consider have no
multiplicative bases for the morphisms spaces in general, see for
instance \cite{bgrs}, a $k$-category is not in general the
linearization of a category. The set theoretical approach cannot
be used, moreover notice that Hochschild-Mitchell (co)homology has
coefficients in a bimodule.
\end{obser}

In case the group $G$ is finite and the characteristic of the
field is zero or do not divide the order of $G$, the algebra $kG$
is semi-simple by Maschke's theorem. Consequently $H_p(G,X)=0$ for
$p>0$ and for any $kG$-module $X$. So at level $2$ the preceding
spectral sequence has differentials zero, which shows the
following.

\begin{prop}
Let ${\cal C} \to {\cal B}$ be a Galois covering of $k$-categories
with finite group $G$ and assume that the characteristic of the
field $k$ is zero or do not divide the order of $G$. Let $\cal M$
be a $\cal B$-bimodule. Then $H_n({\cal B},{\cal M}) = H_n({\cal
C}, L{\cal M})/G$.

Moreover, if $\cal M$ is locally finite, then $H^n({\cal B},{\cal
M}) = H^n({\cal C}, L{\cal M})^G$.
\end{prop}

\begin{obser}

Note that the change of the bimodule of coefficients is
meaningful, namely the lifted bimodule provides the right context.
This is related to a wrong interpretation of the present work made
at the end of the introduction in \cite{marmarmar}, as the example
below  (first considered in \cite{jap-m}) illustrates.

\end{obser}

Consider the free $k$-categories presented by the quivers

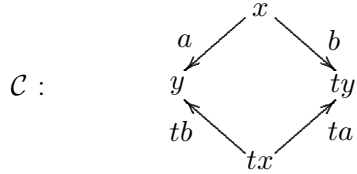
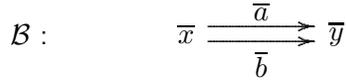
\beginpicture \setcoordinatesystem units <1cm, 1cm> \setplotarea x from -3
to 5, y from -2 to 3.5

\put {$\cal C$ : } [c] at -1 2 \put {$x$} [c] at 2 3
\put {$y$} [c] at .9 2
\put {$ty$} [c] at 3.1 2 \put {$tx$} [c] at 2 1


\arrow <6pt> [.15,.6] from 1.8 2.9 to 1 2.2 \arrow <6pt> [.15,.6]
from 2.2 2.9 to 3 2.2 \arrow <6pt> [.15,.6] from 1.8 1.1 to 1 1.8
\arrow <6pt> [.15,.6] from 2.2 1.1 to 3 1.8

\put {$a$} [r] at 1.1 2.6 \put {$b$} [l] at 2.9 2.6 \put {$tb$}
[r] at 1.1 1.4 \put {$ta$} [l] at 2.9 1.4

\put {$\cal B$ : } [c] at -1 -1 \put {$\overline x$} [c] at 1 -1
\put {$\overline y$} [c] at 3 -1

\arrow <6pt> [.15,.6] from 1.3 -.9 to 2.7 -.9 \arrow <6pt>
[.15,.6] from 1.3 -1.1 to 2.7 -1.1

\put {$\overline a$} [c] at 2 -.7 \put {$\overline b$} [c] at 2
-1.4

\endpicture
\noindent where $\cal B$ is the Kronecker algebra.
There is a Galois covering ${\cal C} \to {\cal B}$ under the
action of the group $G=< t / t^2=1 >$, $\dim H^1({\cal C}, {\cal
C}) =1$ while $\dim H^1({\cal B}, {\cal B}) =3$. But in
characteristic different from two $\dim H^1({\cal C}, L{\cal B})
=5$ and $\dim H^1({\cal C}, L{\cal B})^G =3$ as a simple
computation shows.

We state another immediate consequence of the Cartan-Leray Theorem
for Galois coverings of categories.

\begin{prop}
Let ${\cal C} \to {\cal B}$ be a Galois covering of $k$-categories
with free group $F$ and let $\cal M$ be a $\cal B$-bimodule. Then
the vector spaces $$H_n({\cal B},{\cal M}) \mbox{ and } H_0(F,
H_n({\cal C},L{\cal M})) \oplus H_1(F, H_{n-1}({\cal C},L{\cal
M}))$$ are isomorphic.
\end{prop}

\begin{proof}

Let $F$ be a free group on a set $B$ of generators. A free
resolution of $k$ is given by
$$0 \rightarrow \bigoplus_B kF \rightarrow kF \stackrel{\epsilon}{\rightarrow} k
\rightarrow 0.$$ For a $kF$-module $X$ we have $H_i(kF,X)=0$ for
$i \geq 2$. Consequently the Cartan-Leray sequence has possibly
non zero columns at level $2$ only for $p=0$ and $p=1$, then all
the differentials in this level are zero since they come from $0$
or go to $0$.\qed
\end{proof}

\begin{prop}
With the same assumptions as before, assume in addition that
$\cal M$ is locally finite. Then the vector spaces $$H^n({\cal
B},{\cal M}) \mbox{ and } H^0(F, H^n({\cal C},L{\cal M})) \oplus
H^1(F, H^{n-1}({\cal C},L{\cal M}))$$ are isomorphic.
\end{prop}

\begin{obser}
If $F=T$ is free of rank one the resolution of $k$ is
$$0 \rightarrow kG \stackrel{1-t}{\rightarrow} kG \stackrel{\epsilon}{\rightarrow} k \rightarrow 0$$
Consequently if $X$ is a $kT$-module we have that the cohomology
of $X$ is the kernel and cokernel of the map given by $1-t$, that is,
$H^0(T,X)=X^T$ and $H^1(T,X)=X/T$. The preceding Proposition
specializes as follows for $T$ an infinite cyclic group and $\cal
M$ a locally finite bimodule:

$H^n({\cal B},{\cal M})$ and $H^n({\cal C},L{\cal M})^T \oplus H^{n-1}({\cal C},L{\cal M})/T$ are isomorphic vector
spaces.
\end{obser}

\section {\sf Rank of the Galois group and dimension of $H^1$}

We recall first some well known definitions, see for instance
\cite{ga}. A \emph{basic} $k$-category has no different isomorphic
objects. If $\cal B$ is not basic, choose one object in each
isomorphism class and consider the full subcategory instead. This
corresponds of course to Morita reduction of algebras and
Hochschild-Mitchell (co)homology is invariant under this
procedure, this can be proved by standard arguments.

A \emph{split} $k$-category is a $k$-category $\cal B$ such that
${ }_x{\cal B}_ x$ is a local $k$-algebra for each object $x$.
Notice that if a basic $k$-category is not split but has finite
dimensional endomorphism algebras, the expansion process described
in Section \ref{sectionHOMI} provides a basic split category along
the choice of a complete system of orthogonal primitive
idempotents for each ${ }_x {\cal B}_x$.

In case the field $k$ is algebraically closed and $B$ is a local
finite dimensional $k$-algebra with maximal ideal $\emph{M}$ we
have $B/{\emph M}=k$ and there is a canonical decomposition $B= k
\oplus {\emph M}$. A $k$-category $\cal B$ is \emph{totally split}
if it is split and if for each object $x$ we have such
decomposition ${ }_x {\cal B}_x = k \oplus {\emph M}_x$. If $k$ is
algebraically closed every split category is totally split.

A $k$-category $\cal B$ is \emph{hom-finite} if each vector space
${ }_y {\cal B}_x$ is finite dimensional for each couple of
objects $(y,x)$ of $\cal B$, in other words the standard bimodule
$\cal B$ is locally finite. Notice that we have already quoted
that finite and hom-finite $k$-categories coincide with finite
dimensional $k$-algebras provided with a complete system of
orthogonal idempotents.

Finally a $k$-category is \emph{connected} if the following
equivalence relation $\sim$ has only one class: $\sim$ is
generated by $x \sim y$ if and only if ${}_y{\cal C}_x \not= 0$ or
${}_x{\cal C}_y \not=0$.

We will use the Cartan-Leray spectral sequence of cohomological
type in order to show that for a connected Galois covering with
group $G$ of a totally split, basic and hom-finite $k$-category
$\cal B$, the vector space of group homorphisms $\Hom(G,k)$ is
canonically embedded in $H^1({\cal B},{\cal B})$. As quoted in the
introduction, Assem and de la Pe\~na obtained this result in
\cite{a-jap} when $G$ is the fundamental group of a triangular
finite dimensional algebra presented by a quiver with relations.
In \cite{jap-s} de la Pe\~na and Saor\'\i n noticed that the
triangular hypothesis is superfluous.

We recall first a simple fact from the study of a first quadrant
converging spectral sequence of cohomological type. Namely
$E_{p+1}^{p,0}$ is a subspace of the abutment in degree $p$.
Indeed the differentials at level $p+1$ at the spot (p,0) come from
and go to $0$, which shows that $E_{p+1}^{p,0}=E_{p+2}^{p,0}=
\dots =E_\infty^{p,0}$. In particular $E_2^{1,0}$ is a subspace of
$H^1$.

\begin{teo}
Let ${\cal C} \to {\cal B}$ be a Galois covering of $k$-categories
with group $G$ where $\cal C$ is connected and $\cal B$ is
hom-finite, basic and totally split. Then there is a canonical
vector space monomorphism
$$\Hom(G,k) \hookrightarrow H^1({\cal B}, {\cal B}).$$
\end{teo}

\begin{proof}
Since $\cal B$ is hom-finite we use the Cartan-Leray spectral
sequence of cohomological type converging to $H^n({\cal B}, {\cal
B})$ which has level $2$ term
$$E_2^{p,q}=H^p\left(G, H^q({\cal C}, L{\cal B})\right).$$

At degree one we have a vector spaces exact sequence
$$0\to E^{1,0}_2\to H^1({\cal B}, {\cal B}) \to E^{0,1}_3 \to 0$$
since $E^{0,1}_3=E^{0,1}_4=\cdots E^{0,1}_\infty.$

Notice that $E^{1,0}_2=H^1\left(G, H^0({\cal C}, L{\cal
B})\right)$. We assert that the $kG$-module $H^0({\cal C}, L{\cal
B})$ has a one dimensional $G$-trivial direct summand. Recall that
for any $\cal C$-bimodule $\cal M$ we have
$$H^0({\cal C}, {\cal M}) = \left\{ \left({}_xm_x\right)_{x\in{\cal
C}} \mid {}_yc_x\ {}_xm_x = {}_ym_y\ {}_yc_x \mbox{ for all }
{}_yc_x\in {}_y{\cal C}_x\right\}.$$ Moreover ${}_x(L{\cal B})_x=
{}_{\overline x}{\cal B}_{\overline x}= k\oplus\emph{M}_{\overline
x}$. Consider a family
$$f=(\lambda_x + m_x)_{x\in{\cal C}_0}\in H^0({\cal C}, L{\cal
B}).$$ For $x\neq y$ and in case ${}_y{\cal C}_x$ or ${}_x{\cal C}_y$ is not zero, we have $\lambda_x=\lambda_y$:
indeed let $0\neq c\in {}_y{\cal C}_x$. Since $\cal B$ is basic, ${\overline c}$ is not invertible and it belongs to
the Jacobson radical $r$ of the category $\cal B$ which is is hom-finite. Then there exist a positive integer $i$
such that ${\overline c}\in r^i$ but ${\overline c}\not\in r^{i+1}$. Consider a vector space decomposition
$${ }_{\overline{y}}{\cal B}_{\overline{x}} = k\overline c \oplus
X \oplus\overline{y} r^{i+1} \overline{x}.$$ We have
$$\overline c(\lambda_x+m_x)=\overline c\lambda_x+\overline
cm_x\mbox{ and } (\lambda_y+m_y)\overline c=\lambda_y\overline
c+m_y\overline c$$ and equality holds between these vectors.
Notice that $\overline cm_x\in r^{i+1}$ and $m_y\overline c\in
r^{i+1}$ while $\overline c\lambda_x\in k\overline c$ and
$\lambda_y\overline c\in k\overline c$. Consequently $\overline
c\lambda_x = \lambda_y\overline c$ and $\lambda_x=\lambda_y$.

Since $\cal C$ is connected we infer that $\lambda_x=\lambda_y$
for every pair of objects and $f=(\lambda + m_x)_{x\in{\cal
C}_0}$. Clearly $(\lambda)_{x\in{\cal C}_0}\in H^0({\cal C},
L{\cal B})$, consequently $(m_x)_{x\in{\cal C}_0}\in H^0({\cal C},
L{\cal B})$. We have proved that
$$H^0({\cal C}, L{\cal
B})=k\ \oplus\ \prod_{x\in {\cal C}_0}{\emph M}_{\overline
x}\bigcap H^0({\cal C}, L{\cal B}).$$

Finally we note that this decomposition is a $kG$-modules
decomposition since for each $s\in G$ the corresponding morphism
changes the object spot of the vector space but remains the
identity on the vector space itself. Moreover the resulting action
on $k=(\lambda)_{x \in {\cal C}_0}$ is the identity.

We infer that $H^1(G,k)$ is a canonical direct summand of $E_2^{1,0}$, which in turn is canonically embedded in
$H^1({\cal B},{\cal B})$.

Of course $$H^1(G,k)= \Hom(G, k^+)$$ since derivations with coefficients in a trivial module are just homomorphisms,
and there are no inner derivations besides $0$. \qed
\end{proof}

The exact sequence at the beginning of the preceding proof
provides a description of $H^1(\cal B, B)$ which can be explored
further. In particular conditions for an isomorphism between
$\Hom(G, k^+)$ and $H^1(\cal B, B)$ can be inferred, see
\cite{jap-s}.

\medskip

The following result provides an upper bound for the rank of the
Galois group of a covering, compare with \cite{jap-m} where the
bound is obtained in case of a schurian category $\cal C$ and a
free group acting on $\cal C$ with a finite dimensional basic and
totally split algebra quotient.

\begin{coro}
Let ${\cal C} \to {\cal B}$ be a Galois covering of $k$-categories
with group $G$ where $\cal C$ is connected and $\cal B$ is basic,
hom-finite, and totally split. Then $\rank G \leq \dim
H^1({\cal B},{\cal B})$.
\end{coro}

\begin{proof}
By definition $\rank G = \dimQ (\Q \otimes_\Z
G^{{\mathrm{ab}}})$ where $G^{{\mathrm{ab}}}$ is the
abelianisation of $G$. We have
$$\dimQ (\Q \otimes_\Z G^{{\mathrm{ab}}})= \dimQ \Hom(G,\Q^+)\leq
\dim\Hom (G,k^+) \leq \dim H^1({\cal B}, {\cal B}).$$ \qed
\end{proof}

\footnotesize \noindent C.C.:
\\D\'epartement de Math\'ematiques,
 Universit\'e de Montpellier
2,  \\F--34095 Montpellier cedex 5, France. \\{\tt
Claude.Cibils@math.univ-montp2.fr}

\vskip3mm \noindent M.J.R:\\ Departamento de Matem\'atica,
Universidad Nacional del Sur,\\Av. Alem 1253\\8000 Bah\'\i a
Blanca, Argentina.\\ {\tt mredondo@criba.edu.ar}

\end{document}